\documentclass[onecolumn,11pt,draft=false]{IEEEtran}
\usepackage{mathtools,color}
\usepackage{hyperref}

\usepackage{enumitem}
\usepackage{balance}
\usepackage{setspace}
\usepackage{latexsym}
\usepackage{lipsum}
\usepackage{amsfonts}
\usepackage{amssymb,amsmath,bbm}
\usepackage{cleveref}
\newtheorem{example}{Example}[section]

\newtheorem{lemma}{Lemma}
\newtheorem{theorem}{Theorem}
\newtheorem{proposition}{Proposition}
\newtheorem{corollary}{Corollary}
\newtheorem{proof}{Proof}

\newcommand{\ones}{\mathbf{1}}

\newcommand{\Z}{{\mathbb{Z}}}

\newcommand{\IPCN}{\mathrm{IPCN}}
\renewcommand{\wp}{\textrm{w.p.}}
\newcommand{\oprocendsymbol}{\hbox{$\bullet$}}
\newcommand{\oprocend}{\relax\ifmmode\else\unskip\hfill\fi\oprocendsymbol}
\allowdisplaybreaks

\begin{document}
\title{\LARGE\bf A Finite Memory Interacting P\'{o}lya Contagion Network \\ and its Approximating Dynamical Systems{\LARGE $^*$}}


 \author{Somya Singh,
         Fady Alajaji and
         Bahman Gharesifard}
  \renewcommand{\thefootnote}{\fnsymbol{footnote}}
  \footnotetext{\\ $^*$Department of Mathematics and Statistics, Queen's University, Kingston, Ontario K7L 3N6, Canada (\texttt{17ss262@queensu.ca}, \texttt{fa@queensu.ca}, \texttt{bahman.gharesifard@queensu.ca}). 
  The work of the authors was funded by the Natural Sciences and Engineering Research Council of Canada. The last author wishes to thank 
the Alexander von Humboldt Foundation for their generous support.}
\maketitle




\begin{abstract}
We introduce a new model for contagion spread using a network of interacting \emph{finite memory} two-color  P\'{o}lya urns, which we refer to as the \emph{finite memory interacting P\'{o}lya contagion network}. The urns interact in the sense that the probability of drawing a red ball (which represents an infection state) for a given urn, not only depends on the ratio of red balls in that urn but also on the ratio of red balls in the other urns  in the network, hence accounting for the effect of spatial contagion. The resulting network-wide contagion process is a discrete-time finite-memory ($M$th order) Markov process, whose transition probability matrix is determined. The stochastic properties of the network contagion Markov process are analytically examined, and for homogeneous system parameters, we characterize the limiting state of infection in each urn.
For the non-homogeneous case, given the complexity of the stochastic process, and in the same spirit as the well-studied SIS models, we use a mean-field type approximation to obtain a discrete-time dynamical system for the finite memory interacting P\'{o}lya contagion network. Interestingly, for $M=1$, we obtain a linear dynamical system which exactly represents the corresponding Markov process. For $M>1$, we use mean-field approximation to obtain a nonlinear dynamical system. Furthermore, noting that the latter dynamical system admits a linear variant (realized by retaining its leading linear terms), we study the asymptotic behavior of the linear systems for both memory modes and characterize their equilibrium. Finally, we present simulation studies to assess the quality of the approximation purveyed by the linear and non-linear dynamical systems.
\end{abstract}



\begin{spacing}{1.59}

\section{Introduction}

Interacting urn networks are widely used in the field of applied mathematics, biology, and computer science to model the spread of disease \cite{curing,contagion}, consensus dynamics \cite{consensus,AJ-AM-EM-RS:21}, image segmentation \cite{image} and the propagation of computer viruses in connected devices \cite{internet} and in social networks \cite{social}.
A general description for an interacting two-color urn network is as follows. We are given a network of $N$ urns. At time $t=0$, each urn is composed of some red and some black balls, where different urns can have different initial compositions, but no urn is empty. At each time instant $t$, a ball is chosen for each urn with probability depending on the composition of the urn itself and of the other urns in the network, 
and then additional (``reinforcing'') balls of the color just drawn are added to the urn.
Letting $U_{i,t}$ denote the ratio of red balls in urn $i$ at time $t$, the draw variable $Z_{i,t}$, denoting the indicator function of a red ball chosen for urn $i$ at time $t$, is governed by  
\begin{equation}\label{draw-mechanism}
Z_{i,t}= \begin{dcases}
1 \quad \wp & f(U_{1,t-1},\cdots,U_{N,t-1})\\
0 \quad  \wp & 1-f(U_{1,t-1},\cdots,U_{N,t-1})
\end{dcases}
\end{equation}
where ``\wp''~stands for ``with probability'' and $f: \mathbb{R}^{N}\to (0,1)$ is a real-valued function which accounts for the interaction in the network of urns. The stochastic process $\{Z_{i,t}\}_{t=1}^{\infty}$ is commonly known as a reinforcement process generated by an urn model.

Although, a variety of models are used in the literature depending on the application, P\'{o}lya urns are the most commonly used urn models (there are a few examples based on other interacting urn networks; e.g., see \cite{neeraja,neeraja1} for interacting Friedman urn networks). An interesting feature of P\'{o}lya urns is that their reinforcement scheme represent preferential attachment graphs in the sense that at each time step, we add balls of a particular color with a probability proportional to the number of balls of that color in the urn. This representation of P\'{o}lya urns as preferential attachment graphs makes interacting P\'{o}lya network a good choice to model the spread of infection in an interacting population. For more insight on similarities between P\'{o}lya urns and preferential attachment models see~\cite{internet,collevecchio2013preferential,pa}. 

 In this paper, we study a network of $N$ two-color P\'{o}lya urns. The objective is to model the spread of infection (commonly referred to as contagion) through this interacting P\'{o}lya urn network, where each urn is associated to a node (e.g., ``individual'') in a general network (e.g., ``population") to delineate its ``immunity'' level. Each P\'{o}lya urn in the network contains red and black balls which represent units of ``infection'' and ``healthiness'' respectively. The reinforcement of drawing a ball for each urn is mathematically formulated such that a weighted composition of other urns in the network (activated by an interaction row-stochastic matrix) affects the drawing process, hence capturing interaction between nodes. 

 Various other models have been proposed in the literature to portray contagion in networks using interacting P\'{o}lya processes. In~\cite{contagion}, the concept of ``super urn'' in networks is utilized to model spatial contagion, where drawing a ball for each urn consists of drawing from a ``super urn'' formed by the collection of the urn's own balls and the ones of its neighbours. In \cite{curing,greg}, optimal curing and initialization policies were investigated for the network contagion model of \cite{contagion}. 
 In~\cite{synchronization}, the authors introduce a symmetric type reinforcement scheme for interacting P\'{o}lya urns.
Finally, an example of a more complicated interaction network is given in \cite{ben}, where a finite connected graph with each node equipped with a P\'{o}lya urn is considered and at any given time $t$, only one of the two interacting urns receive balls with probability proportional to a number of balls raised to some fixed power. 
 
 An important characteristic of most reinforcement processes generated via urn  models is that they are non-Markovian in the sense that the composition of each urn at any given time affects its composition at every time instant thereafter. This property is not realistic when modelling the spread of infection as one should account for the possibility that infection is cured (or that the urn is removed, a possibility that we do not consider here). In this work, we consider an interacting P\'{o}lya urn network where each urn has a {\em finite memory}, denoted by $M\geq 1$, in the sense that, at time instant $t>M$, the reinforcing balls added at time $t-M$ are removed from the urn and hence have no effect on future draws. This notion of a finite memory P\'{o}lya urn was introduced in~\cite{finite} (in the context of a single urn) to account for the diminishing effect of past reinforcements on the urn process, which is a realistic assumption when modelling contagion in a population. The resulting network draw variables $\{Z_{t}\}_{t=1}^{\infty} :=\{(Z_{1,t},\cdots,Z_{N,t}\}_{t=1}^{\infty}$ of the P\'{o}lya urn with memory $M$ forms a Markov chain of order $M$. The memory parameter $M$ gives a new degree of freedom to this interacting P\'{o}lya contagion network which makes it more suitable for the study of epidemics. In particular, the draw process of the network with $M=1$ represents the Markov process of the well-known susceptible-infected-susceptible (SIS) model~\cite{hassibi,wang,virus,pare2017epidemic,mei2017dynamics,nowzari2016analysis,andersson2012stochastic}, as each urn $i$ at time $t>1$ exhibits two possible states via its draw variables, susceptible ($Z_{i,t}=0$) and infected ($Z_{i,t}=1$). Hence, our Polya-based model (with $M=1$) is an alternative representation of the SIS model, albeit with different governing parameters. There are other examples in the literature where a Markovian version of the P\'{o}lya process is studied. In~\cite{irene_markov}, a rescaled P\'{o}lya urn model with randomly fluctuating conditional draw probability is considered. Another Markovian P\'{o}lya process is the P\'{o}lya-Lundberg process~\cite{polya-lundberg}, which was recently adapted in~\cite{corona} to measure the dynamics of the SARS-CoV-2 pandemic, among many other models(e.g., see also~\cite{FCF:21} where the classical P\'{o}lya urn scheme is used).
 
 The techniques used in the analysis of a finite memory P\'{o}lya process are quite different from the ones used for any general random reinforcement process. Standard techniques used for the latter case include the method of moments~\cite{friedman}, martingale methods~\cite{synchronization,martingale}, stochastic approximations~\cite{sa,krishanu,sophie} and the embedding of reinforcement processes in continuous time branching processes~\cite{bp,embedding,branching(1)}. A detailed discussion on these methods can be found in the survey~\cite{survey} and in~\cite{irene}. Unlike standard reinforcement processes, we are able to use Markovian properties in our analysis as our model of interacting urns with finite memory $M$ yields an $M$th order Markov draw process. However, one drawback of working with a memory-$M$ Markov chain over a network is that the size of its underlying transition probability matrix grows exponentially with both $M$ and the network size. To account for this problem, after having introduced our interacting P\'{o}lya urn network and investigated its properties in detail, we formulate a {\em dynamical system} to tractably approximate its asymptotic behaviour. To obtain this dynamical system, we make the assumption that for any given time $t>M$, the joint probability distribution of draw processes at times $t-1,\cdots,t-M$ for any urn is equal to the product of marginals. This type of approximation is referred to as ``mean-field approximation'' and is commonly used in the literature on compartmental models, such as the SIS model. The key factor that distinguishes our treatment is the latitude provided by the consideration of memory $M\ge 1$, in contrast to the SIS model which is based on a memory one ($M=1$) Markov chain. In particular, as our simulations which are performed for both non-homogeneous and homogeneous networks verify that the nonlinear dynamical system that we obtain is a good  approximation for the true (underlying) Markov process.
 We also characterize the equilibrium point of this dynamical system, when the nonlinear dynamical system is approximated by its linear part (the latter approximation is exact for the case with memory $M=1$).
 More specifically, we show that when $M=1$,the Markov process mimics precisely a linear dynamical system with a unique equilibrium (which can be exactly determined); while for $M>1$, we note that the approximating linearized dynamical system has a unique equilibrium when its governing (block) matrix has a spectral radius less than unity. 
 In summary, our results provide a novel mathematical framework for the study of epidemics on networks in realistic scenarios where memory is a consideration.
 
 This paper is organised as follows. In Section~\ref{sec:model}, we describe our interacting P\'{o}lya contagion network with memory $M$, which generates a general $M$th order time-varying Markov network draw process. In Section~\ref{sec:homogeneous}, we show that for the homogeneous case (i.e., when all urns have identical initial compositions and the same reinforcement parameters), the network Markov process is time-invariant, irreducible and aperiodic. We obtain the transition probability matrix of this Markov process, illustrate the calculation of its stationary distribution and establish its exact asymptotic marginal distributions. In Section~\ref{sec:three}, we derive the linear dynamical system for the general (non-homogeneous) network for $M=1$ and develop the nonlinear dynamical system approximations for $M>1$ using a mean-field approximation. We investigate the role of memory as well as the equilibrium of these dynamical systems (both for $M=1$ and its linearized variant for $M>1$). Simulation results are presented in Section~\ref{simulations} to assess the modeling quality of the linear and non-linear dynamical systems. Finally, conclusions and directions of future work are stated in Section~\ref{sec:conclusions}.

\section{The Model and its properties}
\label{sec:model} 
We consider a network of $N$ P\'{o}lya urns, where each urn can be associated to a node in an arbitrary network. At time $t=0$, urn $i$ contains $R_{i}$ red balls and $B_{i}$ black balls, $i=1,\ldots,N$. We let $T_{i}= R_{i} + B_{i}$ be the total number of balls in the $i$th urn at time $t=0$, and assume that each urn contains nonzero red and black balls at time $t=0$ i.e., $R_{i}>0$ and $B_{i}>0$. We also let
$U_{i,t}$ denote the ratio of red balls in urn $i$ at time $t$, with its initial value (at time $t=0$) given by $U_{i,0}=R_i/T_i$.

We next define the reinforcement scheme, in the form of {\em draw variables}, $Z_{i,t}$, associated with urn~$i$ at time $t\ge1$, for our proposed interacting P\'{o}lya contagion network:
\[ Z_{i,t}= \begin{cases}
 1 & \textrm{if a red ball is drawn for urn $i$ at time $t$}\\
 0 & \textrm{if a black ball is drawn for urn $i$ at time $t$}
\end{cases}\]
where the {\em process of drawing a ball for urn~$i$} is governed by~\eqref{draw-mechanism} and the function $f$ is explicitly defined below. 
The drawing mechanism~\eqref{draw-mechanism} is applied  {\em simultaneously} to all urns. If a red ball (respectively, a black ball) is drawn for urn $i$, we add $\Delta_{r,i}(t)$ red balls (respectively, $\Delta_{b,i}(t)$ balls) to urn $i$. This scheme, which we refer to as the \emph{urn scheme}, is often captured by a matrix of the form:
$$
    \begin{bmatrix}
    \Delta_{r,i}(t) & 0 \\
    0            & \Delta_{b,i}(t)\\
    \end{bmatrix}.
$$
We assume throughout that $\Delta_{r,i}(t)\geq 0, \Delta_{b,i}(t)\geq 0$, for all $t \in \Z_{\geq 0}$, and that there exist an urn $i$ such that $\Delta_{r,i}(t) + \Delta_{b,i}(t) \neq 0$ at all times $t$.

To formulate the interaction part of the model, we start with defining an \emph{interaction matrix} $S$ as an $N\times N$ row-stochastic matrix with non-negative entries, i.e., each row in $S$ sums to one. Entries of the interaction matrix $S$ are denoted by $s_{ij}$, where~$i,j \in \{1,\ldots,N\}$.
The interaction matrix $S$ can also be thought of as a weighted adjacency matrix of a directed graph with each node equipped with a P\'{o}lya urn. 

Having defined the interaction matrix, we can now explicitly specify the function $f$ used in the drawing mechanism~\eqref{draw-mechanism}. In particular,
we set the probability of choosing a red ball from urn $i$ at time $t$ as follows:
\begin{equation}\label{eq:s-U} 
 Z_{i,t}=\begin{cases}
1 & \wp \quad \sum\limits_{j=1}^{N}s_{ij}U_{j,t-1} \\
0 & \wp \quad 1-\sum\limits_{j=1}^{N}s_{ij}U_{j,t-1}.
\end{cases}
\end{equation}
Note that at time $t>1$, the conditional probability of node~$i$'s draw variable $Z_{i,t}$ is a function of all past draw variables in the network, namely, $\{Z_{j,k}\}$ for $j=1,\ldots,N$ and~$k=1,\ldots,t-1$. Furthermore, as all draws occur simultaneously, the draw variables $Z_{i,t}$ and $Z_{i',t}$ are conditionally independent given all past draws in the network, for any $i\neq i'$; hence at any time $t$,
$$
P\left(Z_{1,t},\ldots,Z_{N,t}|\{Z_{1,k}\}_{k=1}^{t-1},\ldots, \{Z_{N,k}\}_{k=1}^{t-1}\right) \hspace{-0.03in} = \hspace{-0.03in} \prod_{i=1}^N P\left(Z_{i,t}|\{Z_{1,k}\}_{k=1}^{t-1},\ldots, \{Z_{N,k}\}_{k=1}^{t-1}\right).
$$
 For ease of notation, we define the following (normalized) initial and reinforcement network parameters for $i=1,\ldots,N$ and $t\ge1$:
\begin{equation}\label{eq:parameter-IPCN}
    \rho_{i} = \frac{R_{i}}{T_{i}}, \quad \quad \sigma_{i}= 1-\frac{R_{i}}{T_{i}}, \quad \quad \delta_{r,i}(t)=\frac{\Delta_{r,i}(t)}{T_{i}}, \quad \delta_{b,i}(t)=\frac{\Delta_{b,i}(t)}{T_{i}}. 
\end{equation}

We are now ready to describe the notion of {\em finite memory} for the above interacting P\'{o}lya contagion network. In particular, we consider the scenario where we keep the additional balls introduced in the urns after each draw only for a \emph{finite amount of time}, $M\ge 1$, which we call the \emph{memory} of the network. In other words,  the reinforcement process is altered such that, for all urns in the network, we remove the balls added at time $t$ after the $(t+M)$th draw. 
This assumption, which was introduced in~\cite{finite} in the context of a single P\'{o}lya urn, makes the model more realistic, because it accounts for the decrease in severity (or influence) of ``infection'' with time. Also, the balls in the urn at time $t=0$ are never removed from the urn. In the epidemic setting, this initial composition of the urns can be referred to as the \emph{intrinsic or inherent immunity} of the individuals against infection. We will next show that the sequence of $N$-tuple draw variables $\{(Z_{1,t},\ldots,Z_{N,t})\}_{t=1}^{\infty}$ for our finite memory interacting P\'{o}lya network forms an $M$th order Markov chain. For brevity, from now on, we denote our finite memory interacting P\'{o}lya contagion network by $\IPCN(M,N)$, where $M$ is the memory of the network and $N$ is the number of urns in the network; unless we state otherwise, we assume that the underlying parameters are given by~\eqref{eq:parameter-IPCN}. Before showing the Markov property induced by the network-wide draw variable $Z_{t}:= (Z_{1,t},\ldots, Z_{N,t})$, we present a useful characterization of the ratio of red balls at time $t$ in terms of the $\IPCN(M,N)$ parameters in~\eqref{eq:parameter-IPCN}.
\begin{lemma}
\label{lemma:finite_memory}
For an $\IPCN(M,N)$ system, $i=1,2,\ldots,N$, and $t\geq M+1$, we have that
\begin{align}\label{eq:ratio_finite_memory}
U_{i,t}=\frac{\rho_{i} + \sum\limits_{n=t-M+1}^{t}\delta_{r,i}(n)Z_{i,n}}{1 + \sum\limits_{n=t-M+1}^{t}(\delta_{r,i}(n)Z_{i,n} + \delta_{b,i}(n)(1-Z_{i,n}))}.
\end{align}
\end{lemma}
\begin{proof}
Recall that $U_{i,t}$ is the ratio of the red balls in urn $i$ at time $t$. 
Given a finite memory $M$ for the network, at every time instant $t\geq M+1$, we remove the balls added at time $t-M$; hence for $t\geq M+1$, the number of red balls in urn $i$ at time $t$ is
\begin{equation}\label{eq:aux1}
R_{i} + \sum\limits_{n=1}^{t}\Delta_{r,i}(n)Z_{i,n} -\sum\limits_{n=1}^{t-M}\Delta_{r,i}(n)Z_{i,n}
= R_{i} + \sum\limits_{n=t-M+1}^{t}\Delta_{r,i}(n)Z_{i,n}.
\end{equation}
Similarly, the total number of balls in urn $i$ at time $t$ is
\begin{equation}\label{eq:aux2}
T_{i} +  \sum\limits_{n=t-M+1}^{t}(\Delta_{r,i}(n)Z_{i,n} + \Delta_{b,i}(n)(1-Z_{i,n})).
\end{equation}
The result then follows by dividing~\eqref{eq:aux1} by~\eqref{eq:aux2}, and by normalizing both numerator and denominator by $ T_i$.
\end{proof}

We now establish the Markov property for the network draw process  $\{Z_{t}\}_{t=1}^{\infty}$. 
\begin{proposition}\label{lemma:markov_chain}  For an $\IPCN(M,N)$ system, the stochastic process given by $\{Z_{t}\}_{t=1}^{\infty}$ is a time-varying Markov chain of order $M$.
\end{proposition}
\begin{proof}
Let $a_{t}=(a_{1,t},\cdots,a_{N,t})\in \{0,1\}^{N}$. Using~\eqref{eq:s-U} and by virtue of the conditional independence of the draw variables $Z_{i,t}$ and $Z_{i',t}$ given all past draws in the network for all $i\neq i'$, we have for $t\ge M$ that
\begin{align*}
 P[Z_{t+1}=a_{t+1}&|Z_{t}=a_{t},\cdots ,Z_{1}=a_{1}]\\ 
 =& \prod_{i=1}^{N} P[Z_{i,t+1}=a_{i,t+1}|Z_{t}=a_{t},\cdots,Z_{1}=a_{1}]\\ =&\prod\limits_{i=1}^{N}\left(a_{i,t+1}\sum\limits_{j=1}^{N}s_{ij}U_{j,t} + (1-a_{i,t+1})(1-\sum\limits_{j=1}^{N}s_{ij}U_{j,t})\right).
\end{align*}
As a result, we have that 
\begin{align}
\label{eq:tp_general}
P[&Z_{t+1}=a_{t+1}|Z_{t}=a_{t},\cdots ,Z_{1}=a_{1}]\nonumber\\
&=\prod_{i=1}^{N}\Bigg((2a_{i,t+1}-1)\sum\limits_{j=1}^{N}\frac{s_{ij}\Big(\rho_{j} + \sum\limits_{n=t-M+1}^{t}\delta_{r,j}(n)a_{j,n}\Big)}{1 +  \sum\limits_{n=t-M+1}^{t}\left(\delta_{r,j}(n)a_{j,n} + \delta_{b,j}(n)(1-a_{j,n})\right)}  \nonumber \\
& \qquad \qquad \qquad+ (1-a_{i,t+1})\Bigg) \\
&= P[Z_{t+1}=a_{t+1}|Z_{t}=a_{t},\cdots,Z_{t-M+1}=a_{t-M+1}]. \nonumber
\end{align}

 Hence the process $\{Z_{t}\}_{t=1}^{\infty}$ is a time-varying $M$th order Markov chain.
 \end{proof}
 
\section{Analysis of the Homogeneous Case}
\label{sec:homogeneous}
In many settings of contagion propagation, the ``individuals'' in the network, being the urns in our setting , are ``identical'' in the sense of having similar initial parameters. Motivated by this, in this section we develop further the stochastic properties of the $\IPCN(M,N)$ system for the case where the underlying parameters given in~\eqref{eq:parameter-IPCN} are uniform across all urns, by setting
$\Delta_{r,i}(t)=\Delta_{b,i}(t)=\Delta >0, R_{i}= R$ and $T_{i}=T$, for all $i$ and $t$.
 By Lemma \ref{lemma:finite_memory}, for a homogeneous $\IPCN(M,N)$, we have that 
\begin{equation}
\label{eq:homoU}U_{i,t}=\frac{\rho + \delta\sum\limits_{n=t-M+1}^{t}Z_{i,n}}{1+\delta M},
\end{equation}
for $t\ge M+1$, where $\rho = \frac{R}{T}=1-\sigma$ and $\delta = \frac{\Delta}{T}.$

By a proof similar to the one given in Proposition~\ref{lemma:markov_chain}, we next note that for a homogeneous $\IPCN(M,N)$, the draw process $\{Z_{t}\}_{t=1}^{\infty}$ is a {\em time-invariant} Markov chain of order $M$. Indeed, using~\eqref{eq:s-U} and~\eqref{eq:homoU}, or by directly simplifying~ \eqref{eq:tp_general}, we have
\begin{align}
\label{eq:tp}
\mathbf{P}[Z_{t+1}=a_{t+1}&|Z_{t}=a_{t},\cdots,Z_{1}=a_{1}]\nonumber\\
&=\prod_{i=1}^{N}\left((2a_{i,t+1}-1)\sum_{j=1}^{N}\frac{s_{ij}(\rho + \delta \sum_{k=t-M+1}^{t}a_{j,k})}{1+\delta M} + (1-a_{i,t+1})\right)\nonumber\\
&= \mathbf{P}[Z_{t+1}=a_{t+1}|Z_{t}=a_{t},\cdots,Z_{t-M+1}=a_{t-M+1}],
\end{align}
where the conditional probabilities do not depend on time; hence for the homogeneous $\IPCN(M,N)$, $\{Z_{t}\}_{t=1}^{\infty}$ is a time-invariant $M$th order Markov chain.
 
We next examine more closely the properties of this time-invariant Markov chain. Starting with the case of $M=1$, i.e., for the homogeneous $\IPCN(1,N)$ system, the transition probabilities of $\{Z_{t}\}_{t=1}^{\infty}$ are given by~\eqref{eq:tp}. Also, the probability of going from state $a = (a_{1,t},\cdots,a_{N,t})$ to state $b=(b_{1,t+1},\cdots ,b_{N,t+1})$ can expressed as
\[q_{ab}^{(1,N)} := q_{ab}^{(1)}q_{ab}^{(2)}\cdots q_{ab}^{(N)}\]
where
 \begin{equation}
     \label{eq:tpM=1}q_{ab}^{(d)}=\begin{dcases} 
      \frac{\sigma + (1-\sum\limits_{k=1}^{N}s_{dk}a_{k,t})\delta}{1+ \delta} & \textrm{if} \hspace{5 mm} b_{d,t+1}=0 \\
      \frac{\rho + \sum\limits_{k=1}^{N}s_{dk}a_{k,t}\delta}{1+ \delta} & \textrm{if}\hspace{5 mm}   b_{d,t+1}= 1,
      \end{dcases}
      \end{equation}
      with $ d\in \{1,\cdots, N\}$. We denote the transition probability matrix of this Markov process with memory $M=1$ and $N$ urns by the matrix $Q^{(1,N)}=\big[q_{ab}^{(1,N)}\big]$, whose entries are given above.
Note that for a memory $M=1$ Markov process, it is possible to go from any state to any state in one time step with a positive transition probability. Hence, the Markov chain is irreducible and aperiodic.

For the case of $M>1$, since $\{Z_{t}\}_{t=1}^{\infty}$ is an $M$th order Markov process with~$2^N$ states, the process $\{W_t\}_{t=1}^{\infty}$, defined by $W_{t}:=(Z_{t},Z_{t+1},\cdots ,Z_{t+M-1})$, becomes an (equivalent) Markov chain of order one with an expanded alphabet of $2^{NM}$ states. For the Markov chain $\{W_t\}$,
the transition probability of going from state 
\[
a = ((a_{11},a_{21},\cdots ,a_{N1}),\cdots, (a_{1M},a_{2M},\cdots ,a_{NM}))
\]
to  state 
\[
b = ((b_{11},b_{21}
,\cdots ,b_{N1}),\cdots ,(b_{1M},b_{2M},\cdots ,b_{NM}))
\]
in one time step is non zero if and only if  $a_{ij}=b_{i(j-1)}$ for  $i \in \{1,\cdots,N\}$ and $j \in \{2,\cdots ,M\}$. If the transition probability is nonzero, it is given by
 \[
      q^{(M,N)}_{ab}:=
  \widetilde{q}_{ab}^{(1)}\widetilde{q}_{ab}^{(2)}\cdots \widetilde{q}_{ab}^{(N)}
 \]
 where
  \begin{equation}\label{eq:tp>1}\widetilde{q}_{ab}^{(d)}=\begin{dcases} 
      \frac{\sigma + \Big(M-\sum\limits_{i=1}^{N}s_{di}(\sum\limits_{k=1}^{M}a_{ik})\Big)\delta}{1+ M\delta} & \textrm{if} \hspace{5 mm} b_{dM}=0 \\\\
      \frac{\rho + \Big(\sum\limits_{i=1}^{N}s_{di}(\sum\limits_{k=1}^{M}a_{ik})\Big)\delta}{1+ M\delta} & \textrm{if}\hspace{5 mm}   b_{dM}= 1,
      \end{dcases} 
\end{equation}
where $ d\in \{1,\cdots, N\}$.
 Similar to the case of memory one, we denote the transition probability matrix for memory $M$ and $N$ urns by the matrix $Q^{(M,N)}=\big[q_{ab}^{(M,N)}\big]$ whose entries are given above. In the next lemma, we extend the above irreducibility and aperiodicity properties for the Markov process $\{W_{t}\}_{t=1}^{\infty}$ with $M\geq 1$.
 \begin{lemma}\label{lemma:irreducible}
 For the homogeneous $\IPCN(M,N)$, the transition probability matrix $ Q^{(M,N)}$ is irreducible and aperiodic.
 \end{lemma}
 \begin{proof}
 We have already seen that the Markov process given by the transition probability matrix $Q^{(1,N)}$ is irreducible and aperiodic. For memory $M>1$, to prove irreducibility of the Markov chain, we show that given any two states, it is possible to go from one state to another in finitely many time steps with a positive probability. Let us fix two arbitrary states, $ a = ((a_{11},a_{21},\cdots ,a_{N1}),\cdots ,(a_{1M},a_{2M},\cdots ,a_{NM})) $ and
 $ b = ((b_{11},b_{21},\cdots,b_{N1}),\cdots,(b_{1M},b_{2M},\cdots,b_{NM})) $. We next construct an $M$-step path (which occurs with a positive probability) between states $a$ and $b$.
    \begin{itemize}
        \item Suppose the  Markov chain is in state $W_{t}=a$ at time $t$. At time $t+1$, we go from state $a$ to state,
        \begin{align*}
        W_{t+1}&=a^{(0)}=\Big( Z_{t+1}=(a_{12}, a_{22},\cdots, a_{N2}),\\
        &\cdots,Z_{t+M}=(a_{1M},a_{2M},\cdots,a_{NM}),
         Z_{t+M}=(b_{11},b_{21},\cdots ,b_{N1})\Big).
         \end{align*} 
        Since $a_{ij}=a^{(0)}_{i(j-1)}$ for $i \in \{1,2,\cdots,N\}$ and $j \in \{2,3,\cdots,M\}$, the transition probability of going from state $a$ to $a^{(0)}$ is nonzero and can be obtained using  \eqref{eq:tp>1}. 
        \item At time $t+2$ we go from state $a^{(0)}$ to state $a^{(1)}$
         \begin{align*}W_{t+2}=a^{(1)}&=\Big(Z_{t+2}=(a_{13},\cdots ,a_{N3}),\cdots,Z_{t+M+1}=(b_{11},\cdots,b_{N1}),\\
        &Z_{t+M+1}=(b_{12},\cdots,b_{N2})\Big).\end{align*}
    \end{itemize}
     Following this pattern of adding one $N$-tuple from state $b$ at each time step, we will reach state $b$ in $M$ time steps.
     In summary, choosing any initial state, we can reach any other state of the Markov chain in at most $M$ steps. Hence, the Markov chain is irreducible. Also, note that the period of the state with all zeros is one. Since all the states of an irreducible Markov chain have the same period, we obtain that this Markov chain is aperiodic. 
\end{proof} 
 
 We now have an explicit formula for the entries of the transition probability matrix $Q^{(M,N)}$. Since the time-invariant Markov chain $\{W_{t}\}_{t=1}^{\infty}$ is irreducible and aperiodic, it has a unique stationary distribution and it is ergodic. We next illustrate this Markov process via a simple example. 
 \begin{example}\label{example:stationarydistribution1}
 \em{
 Given a homogeneous $\IPCN(1,2)$ system with interaction matrix
 \[S=\begin{bmatrix}
 s_{11} & 1-s_{11}\\
 s_{21}  & 1-s_{21}
 \end{bmatrix},\]
 the stationary distribution for the transition probability matrix $Q^{(1,2)}$ is given by
 \begin{align*}
\pi_{00}&=\frac{2\sigma^{2}\delta + \sigma^{2} + (1-s_{11}-s_{21}+2s_{11}s_{21})\sigma\delta^{2}}{(1-s_{11}-s_{21}+2s_{11}s_{21})\delta^{2}+2\delta+1}\\
\pi_{01}&=\frac{\rho\sigma(1+2\delta)}
{(1-s_{11}-s_{21}+2s_{11}s_{21})\delta^{2}+ 2\delta + 1}\\
 \pi_{10}&=\frac{\rho\sigma(1+2\delta)}{(1-s_{11}-s_{21}+2s_{11}s_{21})\delta^{2}+2\delta+1}\\
\pi_{11}&=\frac{\rho(2\delta-\sigma-2\sigma\delta+ (1-s_{11}-s_{21}+2s_{11}s_{21})\delta^{2}+1)}{(1-s_{11}-s_{21}+2s_{11}s_{21})\delta^{2}+2\delta+1}.
\end{align*}
It is easy to see that $\pi = [\pi_{00}, \pi_{01},\pi_{10},\pi_{11}] $  
satisfies the equation $ \pi Q^{(1,2)}=\pi$.
We thus have
\begin{itemize}\item $\lim\limits_{t \to \infty}P(Z_{1,t}=1)=\pi_{10}+\pi_{11}=\rho,  \lim\limits_{t \to \infty}P(Z_{2,t}=1)=\pi_{01} + \pi_{11}=\rho$,
    \item $\lim\limits_{t \to \infty}P(Z_{1,t}=0)=\pi_{00} + \pi_{01}=\sigma,  \lim\limits_{t \to \infty}P(Z_{2,t}=0)=\pi_{00} + \pi_{10}=\sigma$.
    \end{itemize}
Also using \eqref{eq:homoU} with $M=1$, we have for $i=1,2$ that
\[\lim_{t \to \infty}\mathbb{E}[U_{i,t}] = \frac{\rho + \delta \lim\limits_{t\to \infty}\mathbb{E}[Z_{i,t}]}{1+ \delta } 
=\frac{\rho + \delta \lim\limits_{t\to \infty}P(Z_{i,t}=1)}{1+ \delta }=
\rho.\]

Hence, irrespective of the used interaction matrix, the asymptotic marginal (one-fold) distributions and urn compositions for the $\IPCN(1,2)$ system are the same as for the single (memory one) P\'{o}lya urn studied in~\cite{finite}; this result is proved in general in Theorem~\ref{thm:homoipcn} below.
We, however, next observe that the asymptotic $2$-fold draw distributions for the $\IPCN(1,2)$ urns do not match their counterparts for the single P\'{o}lya urn process of~\cite{finite}. 
Indeed, the 2-fold (joint) distribution vector of the single-urn (stationary) P\'{o}lya Markov chain in~\cite{finite} is given~by
 \[\tilde{\pi}^{(2)}=\left[\frac{\sigma(\sigma+\delta)}{1+\delta},  \frac{\rho\sigma}{1+\delta}, \frac{\rho\sigma}{1+\delta}, \frac{\rho(\rho+\delta)}{1+\delta}\right].\]
Furthermore, for the homogeneous $\IPCN(1,2)$ system,
the joint probability $P(Z_{1,t}=a_{1},Z_{1,t+1}=b_{1})$ for urn~1 of the homogeneous $\IPCN(1,2)$ system is given by
 \begin{align*}
 &P(Z_{1,t}=a_{1},Z_{1,t+1}=b_{1})\\
 &=\sum\limits_{a_{2},b_{2}\in\{0,1\}}P(Z_{1,t}=a_{1},Z_{1,t+1}=b_{1},Z_{2,t}=a_{2},Z_{2,t+1}=b_{2})\\
 &=\sum\limits_{a_{2},b_{2}\in\{0,1\}}P(Z_{1,t+1}=b_{1},Z_{2,t+1}=b_{2}|Z_{1,t}=a_{1},Z_{2,t}=a_{2}) P(Z_{1,t}=a_1,Z_{2,t}=a_2).
 \end{align*}
 Thus noting the conditional independence of $Z_{1,t+1}$ and
 $Z_{2,t+1}$ given $(Z_{1,t},Z_{2,t})$, and using the $\IPCN(1,2)$ matrix $Q^{(1,2)}$
 along with the fact that 
 \[
 \lim_{t \to \infty}P(Z_{1,t}=a_1,Z_{2,t}=a_2)=\pi_{a_{1},a_{2}},
 \]
 we obtain
\begin{align*}
  \lim_{t \to \infty}P(Z_{1,t}=0,Z_{1,t+1}=0)&=\frac{\sigma(\sigma+\delta)}{1+\delta}-\frac{\pi_{01}(1-s_{11})\delta}{(1+\delta)}\\
\lim_{t \to \infty}P(Z_{1,t}=0,Z_{1,t+1}=1)&=\frac{\sigma\rho}{1+\delta}
+\frac{\pi_{01}(1-s_{11})\delta}{(1+\delta)}\\
 \lim_{t \to \infty}P(Z_{1,t}=1,Z_{1,t+1}=0)&=\frac{\sigma\rho}{1+\delta}
 +\frac{\pi_{01}(1-s_{11})\delta}{(1+\delta)}\\
 \lim_{t \to \infty}P(Z_{1,t}=1,Z_{1,t+1}=1)&=\frac{\rho(\rho+\delta)}{1+\delta} 
 -\frac{\pi_{01}(1-s_{11})\delta}{(1+\delta)},
 \end{align*}
 which shows explicitly by how much the asymptotic $2$-fold draw distribution for urn~1 deviates from $\tilde{\pi}^{(2)}$. Note that by setting $s_{11}=1$, the error term $\pi_{01}(1-s_{11})\delta/(1+\delta)$ reduces to zero,
 making the two distributions match, as expected (since when
 $s_{11}=1$, urn~1 only interacts with itself).

}
\oprocend
\end{example}
Note that it is much harder to derive in closed-form the stationary distribution for the homogeneous $\IPCN(M,N)$ system with $M > 1 $ and $N >2$ but we have the following asymptotic marginal probabilities for a homogeneous $\IPCN(M,N)$ system.
\begin{theorem}\label{thm:homoipcn}
For a homogeneous $\IPCN(M,N)$ system 
\begin{equation}\label{eq:marginal}
\lim_{t \to \infty}P(Z_{i,t}=1)=\rho
\end{equation}
for all urns $i$ in the network.
\end{theorem}
\begin{proof}
Let $\gamma_{i}=\lim\limits_{t \to \infty}E[Z_{i,t}]$ for $i \in \{1,2,\cdots,N\}$. Using \eqref{eq:tp>1}, we obtain
\[\gamma_{d} = \frac{\rho + \sum_{i=1}^{N}s_{di}(M\gamma_{i})\delta}{1 + M\delta}\]
for $d\in \{1,2,\ldots,N\}$.
Let $\textbf{1}_{N}=[1,\cdots,1]^{\textbf{T}}$ and $\gamma = [\gamma_{1},\cdots,\gamma_{N}]^{\textbf{T}}$, where $^{\textbf{T}}$ denotes transposition. Then,
\[(1 + M\delta)\gamma = \rho \textbf{1}_{N} + (M\delta)S\gamma\]
which gives
\[(1 + M\delta)(\gamma - \rho\textbf{1}_{N}) = (M\delta)S(\gamma - \rho\textbf{1}_{N}).\]
Setting $\widetilde{\gamma} :=\gamma -\rho \textbf{1}_{N}$ in the above equation, we have that
\[S\widetilde{\gamma} = \frac{1+M\delta}{M\delta}\widetilde{\gamma}.\]
Since the eigenvalues of $S$ have absolute values less than or equal to one (as $S$ is a row-stochastic matrix), we obtain that
\[\widetilde{\gamma}=0\]
which implies that
\[\gamma_{i} =\rho \quad \forall \quad i \in  \{1,2\cdots,N\}.\]
\end{proof}
\section{Dynamical System Models}
\label{sec:three}
As seen in the earlier section, in general it is not easy to obtain the stationary distribution for the $\IPCN(M,N)$ Markov chain characterized in \eqref{eq:tp_general}, which has $2^{MN}$ states. Due to this exponential increase in the size of the transition probability matrix with the number of urns $N$ and memory $M$, it is difficult to analytically solve for the stationary distribution in terms of the system parameters. In this section, we present the main core of our paper, i.e., a class of dynamical systems whose trajectory approximates the infection probability at time $t$ for any urn $i$ in the network. Notably, for the case of $M=1$, we observe that the process naturally leads to an exact linear dynamical system without using any approximation. For $M>1$, we use a \emph{mean-field approximation} to obtain an approximating dynamical system. This is in the same theme as the classical SIS model, and its variations, where a dynamical system is often used instead of the original Markov chain. A few examples in which an approximate dynamical system is constructed for a Markov chain are given in \cite{hassibi, wang, virus,andersson2012stochastic}. Throughout the section, we consider $\delta_{r,i}(t)=\delta_{r,i}$ and $\delta_{b,i}(t)=\delta_{b,i}$ for all time instances $t$, i.e., we remove the time dependence from the reinforcement parameters of the Markov process.

\subsection{Dynamical system for $\textbf{M=1}$}
 Our main objective is to analyze the behaviour of the draw variables $\{Z_{i,t}\}$, when memory is one. In particular, we obtain a dynamical system for the evolution in time of $P(Z_{i,t}=1)$, which we outline next. 
 
 For ease of notation, given an $\IPCN(M,N)$ and an urn $i$, we denote the infection probability at time $t$ by
\[ P_{i}(t):=P(Z_{i,t}=1).\]
 Recall from Lemma~\ref{lemma:finite_memory} and \eqref{eq:tp_general} with $M=1$ that the conditional infection probability of urn $i$ at time $t$, given all the draw variables at time $t-1$, is given by
 \begin{align}\label{eq:M=1}
 P(Z_{i,t}=1| Z_{1,t-1},Z_{2,t-1},\cdots ,Z_{N,t-1}) &= \sum_{j=1}^{N}\frac{s_{ij}(\rho_{j} + \delta_{r,j}Z_{j,t-1})}{1 + \delta_{r,j}Z_{j,t-1} + (1-Z_{j,t-1})\delta_{b,j}}\\
 &= \sum_{j=1}^{N}[s_{ij}\beta_{1}^{(j)}(1)Z_{j,t-1} + s_{ij}\beta_{1}^{(j)}(0)(1-Z_{j,t-1})]\nonumber
 \end{align}
where 
\[\beta^{(j)}_{1}(k) := \dfrac{\rho_{j} + k\delta_{r,j} }{1 + k\delta_{r,j} + (1-k)\delta_{b,j}}, \quad j \in \{1,\cdots,N\},\quad k \in \{0,1\}.\]
Now taking expectation with respect to $(Z_{1,t-1},\cdots,Z_{N,t-1})$ on both sides of \eqref{eq:M=1}, we get
\begin{align}\label{eq:M=1single}
P_{i}(t) = \sum_{j=1}^{N}[\beta_{1}^{(j)}(1)s_{ij}P_{j}(t-1) + s_{ij}\beta_{1}^{(j)}(0)(1-P_{j}(t-1))].
\end{align}
To this end, defining the vector $P(t)$ as
\[P(t)= [
P_{1}(t),
P_{2}(t),\cdots,
P_{N}(t)
]^{\mathbf{T}},\]
we obtain the following dynamical system for the $\IPCN(1,N)$ network.
\begin{theorem}\label{thm:memory1}
 For the $\IPCN(1,N)$ system, the infection vector satisfies the equation
\begin{equation}\label{eq:linear_dynamical_system}P(t) =  J_{N,1}P(t-1) + C_{N,1}
\end{equation}
where $J_{N,1} \in \mathbb{R}^{N \times N}$, $C_{N,1} \in \mathbb{R}^{N \times 1}$ are matrices with respective entries:
\[[J_{N,1}]_{i\times j}= \frac{s_{ij}(\rho_{j}+\delta_{r,j})}{(1+\delta_{r,j})}-\frac{s_{ij}\rho_{j}}{(1+\delta_{b,j})} = s_{ij}(\beta^{(j)}_{1}(1)-\beta^{(j)}_{1}(0))\]
\[\textrm{and} \quad [C_{N,1}]_{1\times i}=\sum\limits_{j=1}^{N}\frac{s_{ij}\rho_{j}}{(1+\delta_{b,j})} = \sum\limits_{j=1}^{N}s_{ij}\beta^{(j)}_{1}(0).\]
\end{theorem}

\begin{proof}
Follows from \eqref{eq:M=1single}.
\end{proof}

  We next examine the equilibrium of this linear dynamical system.
\begin{theorem}\label{thm:stability}
The linear dynamical system for the $\IPCN(1,N)$ system given by~\eqref{eq:linear_dynamical_system}
has a unique equilibrium point given by $P^{*}=(I-J_{N,1})^{-1}C_{N,1}$ and 
\[\lim\limits_{t\to \infty} P_{i}(t)= P^{*}_i\]
for all $ i\in\{1,\ldots,N\}$.
\end{theorem}
\begin{proof}
It is enough to show that the spectral radius of the matrix $J_{N,1}$ is less than one; since the spectral radius is less than, or equal to, the row sum norm of the matrix, it is enough to show that the row sum norm of $J_{N,1}$ is strictly less than $1$, see~\cite{spectrum}. 
Note that
\begin{equation}\label{eq:stability}
-1 < \frac{(\rho_{j}+ \delta_{r,j})}{(1+\delta_{r,j})}- \frac{\rho_{j}}{(1+\delta_{b,j})} < 1 
\end{equation}
since $0\leq \rho_{j}\leq 1$ and   $\delta_{r,j},\delta_{b,j}\geq 0$ for all $j \in \{1,2...,N\}.$
Hence, the sum of absolute values of entries in $i$th row of the matrix $J_{N,1}$ satisfies
\[\sum_{j=1}^{N}s_{ij}\left|\frac{(\rho_{j}+ \delta_{r,j})}{(1+\delta_{r,j})}- \frac{\rho_{j}}{(1+\delta_{b,j})}\right| < \sum_{j=1}^{N}s_{ij} = 1,\]
which yields the result.
\end{proof}
As an illustration, we find the  equilibrium of the linear dynamical system \eqref{eq:linear_dynamical_system} for a much simpler $\IPCN(1,N)$ system.

\begin{corollary}\label{corollary}
Given an $\IPCN(1,N)$ system with $S=I$,
\begin{equation}\label{eq:fixed_point_ind}\lim_{t \to \infty}P_{i}(t)= \frac{\rho_{i}(1+\delta_{r,i})}{1+\delta_{b,i} + \rho_{i}(\delta_{r,i}-\delta_{b,i})}.\end{equation}

\end{corollary}

\begin{proof}
 For an $\IPCN(1,N)$ system with $S=I$, we have that
\[P(Z_{1,t-1}=a_{1},\cdots ,Z_{N,t-1}=a_{N}) =  \prod\limits_{j=1}^{N}P(Z_{j,t-1}=a_{j}). \]
In this case, since $S=I$ and hence the draw variables of urns are independent of each other.
The asymptotic value of $P_{i}(t)$ for $ i \in \{1,\cdots,N\}$ is given by the equilibrium point of the linear dynamical system
\[P(t) = J_{N,1}P(t-1) + C_{N,1}\]
 which is given by $P^{*}\in \mathbb{R}^{N}$ whose $i$th component is given by \eqref{eq:fixed_point_ind}.
 
Another way to find this equilibrium point is to write the transition probability matrix for a single urn using \eqref{eq:tp_general} and solving for stationary distribution to obtain $\lim\limits_{t \to \infty}P_{i}(t)$. The transition probability matrix for a single non-homogeneous urn $i$ is given by
\[Q^{(1,1)}=\begin{bmatrix}
\dfrac{\sigma_{i}+\delta_{b,i}}{1+\delta_{b,i}}& \dfrac{\rho_{i}}{1+\delta_{b,i}}\\
\dfrac{\sigma_{i}}{1+\delta_{r,i}}& \dfrac{\rho_{i}+\delta_{r,i}}{1+\delta_{r,i}}
\end{bmatrix}.\]
On solving for the stationary distribution, $[\pi_{0}, \pi_{1}]Q^{(1,1)}=[\pi_{0}, \pi_{1}]$, we obtain that $\pi_{1}$ indeed equals the right-hand-side (R.H.S.) of~\eqref{eq:fixed_point_ind}.
\end{proof}

We also illustrate Theorem~\ref{thm:stability} by examining the special homogeneous case. This aligns with the result in Theorem~\ref{thm:homoipcn}.
\begin{corollary}
\label{fixedpoint}
For a homogeneous $\IPCN(1,N)$ system, the  equilibrium  of \eqref{eq:linear_dynamical_system}
 is given by $P^{*}= \rho \ones_N$, where $ \ones_N$ is vector of ones of size $N$.
 \end{corollary}
\begin{proof}
By Theorem~\ref{thm:stability}, the equilibrium $P^{*}$ is given by
\[ P^{*} = (I-J_{N,1})^{-1}C_{N,1}.\]
Note that the row sums in $(I-J_{N,1})$ are given by $\dfrac{1}{1+\delta}$, i.e.,
\[(I-J_{N,1})\ones_N = \frac{1}{1+\delta}\ones_N.\]
Therefore,
\[(I-J_{N,1})^{-1}C_{N,1} = (I-J_{N,1})^{-1}
\begin{bmatrix}
\frac{\rho}{(1+\delta)}&
\frac{\rho}{(1+\delta)}&
\cdots &
\frac{\rho}{(1+\delta)}
\end{bmatrix}^{\mathbf{T}} = \rho \ones_N.\]

\end{proof}

\subsection{Dynamical system for $\mathbf{M>1}$}
We now construct a class of dynamical systems which approximates the $\IPCN(M,N)$ Markov chain in \eqref{eq:tp_general}. Unlike the memory $M=1$ case, we need to resort to approximations here to obtain dynamical systems. We use the following \emph{mean-field approximation here:} 
   \begin{itemize}
    \item We assume that for every time instant~$t>M$, for each urn $i$, $Z_{i,t-1}$,$\cdots$,$Z_{i,t-M}$ are approximately independent of each other; i.e., at any given time instant~$t>M$, we assume that
\begin{align}\label{eq:mean_field}
P[Z_{j,t-1},Z_{j,t-2},\cdots,Z_{j,t-M}] \approx \prod_{k=1}^{M}P[Z_{j,t-k}],
\end{align}
for all $j \in \{1,2,\cdots,N\}$.
\end{itemize}

\medskip
For the $\IPCN(M,N)$ system, we have from \eqref{eq:tp_general} that
\begin{align}\label{eq:M>1}P[Z_{i,t}=&1|(Z_{1,t-1},\cdots ,Z_{1,t-M}),\cdots ,(Z_{N,t-1},\cdots,Z_{N,t-M})] = \nonumber\\& \sum\limits_{j=1}^{N}\frac{s_{ij}(\rho_{j} + \delta_{r,j}\sum\limits_{k=1}^{M}Z_{j,t-k})}{1+\sum\limits_{k=1}^{M}(\delta_{r,j}Z_{j,t-k} + \delta_{b,j}(1-Z_{j,t-k}))}.
\end{align}
Now, taking expectation with respect to $$((Z_{1,t-1},\ldots,Z_{1,t-M}), \ldots,((Z_{N,t-1},\ldots,Z_{N,t-M}))$$
on both sides of \eqref{eq:M>1} and using the linearity property of expectation, we obtain
\begin{align}\label{eq:dymanicalsystem}
&P[Z_{i,t}=1] = \sum_{j=1}^{N}E\left[\frac{s_{ij}(\rho_{j} + \delta_{r,j}\sum\limits_{k=1}^{M}Z_{j,t-k})}{1 + \sum\limits_{k=1}^{M}(\delta_{r,j}Z_{j,t-k} + \delta_{b,j}(1-Z_{j,t-k}))}\right]\\ \nonumber
& = \sum_{j=1}^{N}\sum_{B_{M}}\frac{s_{ij}(\rho_{j} + \delta_{r,j}\sum_{k=1}^{M}a_{k})}{1+ \sum\limits_{k=1}^{M}(\delta_{r,j}a_{k} + \delta_{b,j}(1-a_{k}))}P(Z_{j,t-1}=a_{1},\cdots Z_{j,t-M}=a_{M})
\end{align}
where
\begin{align}\label{eq:notation}
B_{M} := \{(a_{1},a_{2},\cdots a_{M}) \ | \ 
  a_{k}\in \{0,1\} \quad  \textrm{for} \quad  k \in \{1,2,\cdots ,M\}\}.
\end{align}


Now we use the \emph{mean-field approximation} \eqref{eq:mean_field} in \eqref{eq:dymanicalsystem} to obtain the following class of approximating nonlinear dynamical systems
\begin{align}\label{eq:M>1(2)}
&P_{i}(t) \approx \sum_{j=1}^{N}\sum_{B_{M}}\frac{s_{ij}(\rho_{j} + \delta_{r,j}\sum_{k=1}^{M}a_{k})}{1+ \sum\limits_{k=1}^{M}(\delta_{r,j}a_{k} + \delta_{b,j}(1-a_{k}))}\prod_{k=1}^{M}P[Z_{j,t-k}=a_{k}]=\\ \nonumber
& \sum_{j=1}^{N}\sum_{B_{M}}\frac{s_{ij}(\rho_{j} + \delta_{r,j}\sum_{k=1}^{M}a_{k})}{1+ \delta_{r,j}\sum\limits_{k=1}^{M}a_{k} + \delta_{b,j}(M-\sum\limits_{k=1}^{M}a_{k})}\prod_{k=1}^{M}(a_{k}P_{j}(t-k) + (1-a_{k})(1-P_{j}(t-k))).
\end{align}
For simplicity of notation, we write \eqref{eq:M>1(2)} in the following way:
\begin{align}\label{eq:nonlineards}
P_{i}(t) \approx \sum_{B_{M}}[\sum_{j=1}^{N}s_{ij}\beta_{M}^{(j)}(v_{M})]\prod_{k=1}^{M}(a_{k}P_{j}(t-k) + (1-a_{k})(1-P_{j}(t-k)))
\end{align}
where
\begin{align}\label{eq:notation}
&v_{M}=\sum\limits_{k=1}^{M}a_{k} \quad \textrm{where} \quad a_{k}\in\{0,1\} \quad \textrm{for} \quad k \in \{1,2,\cdots ,M\},\nonumber\\
& \beta^{(j)}_{M}(l) = \dfrac{\rho_{j} + l\delta_{r,j} }{1 + l\delta_{r,j} + (M-l)\delta_{b,j}}, \quad j \in \{1,\cdots,N\},l \in \{0,1,\cdots M\}.
\end{align}
We next give a useful rearranged form of \eqref{eq:nonlineards}. In particular, even though~\eqref{eq:nonlineards} appears to be complicated, after some simplifications, the coefficients of the nonlinear terms follow a binomial pattern. To give an idea of this binomial pattern, we will first present a few examples and then give a proof formula for the rearranged form of \eqref{eq:nonlineards}. 

\begin{example}
\em{
We note that for the $\IPCN(2,2)$ system, the approximating dynamical system is given by
\begin{align*}
P_{i}(t) \approx {}& 
\sum\limits_{j=1}^{2}s_{ij}\beta^{(j)}_{2}(0)
 + \sum_{j=1}^{2}\sum_{k=1}^{2}s_{ij}\left(\beta^{(j)}_{2}(1)-\beta^{(j)}_{2}(0)\right)P_{j}(t-k) \\
& + \sum\limits_{j=1}^{2}s_{ij}\left(\beta^{(j)}_{2}(2)-2\beta^{(j)}_{2}(1)+\beta^{(j)}_{2}(0)\right)\prod\limits_{k=1}^{2}P_{j}(t-k).
\end{align*}
Next, by expansion, we observe that for $\IPCN(3,2)$ system, the approximating dynamical system is given by
\begin{align*}
    P_{i}(t)&\approx{}
    \sum\limits_{j=1}^{2}s_{ij}\beta^{(j)}_{3}(0) + \sum\limits_{j=1}^{2}\sum\limits_{k=1}^{3}P_{j}(t-k)s_{ij}\left[\beta^{(j)}_{3}(1)-\beta^{(j)}_{3}(0)\right] \\
    &+ \sum\limits_{j=1}^{2}P_{j}(t-1)P_{j}(t-2)s_{ij}\left[\beta^{j}_{3}(0)-2\beta^{(j)}_{3}(1)+\beta^{(j)}_{3}(2)\right]\\
& + \sum\limits_{j=1}^{2}P_{j}(t-2)P_{j}(t-3)s_{ij}\left[\beta^{(j)}_{3}(0) -2\beta^{(j)}_{3}(1)+\beta^{j)}_{3}(2)\right]\\
& + \sum\limits_{j=1}^{2}P_{j}(t-1)P_{j}(t-3)s_{ij}\left[\beta^{(j)}_{3}(0)-2\beta^{(j)}_{3}(1)+\beta^{(j)}_{3}(2)\right]\\
& + \sum\limits_{j=1}^{2}P_{j}(t-1)P_{j}(t-2)P_{j}(t-3)s_{ij}\left[3\beta^{(j)}_{3}(1)-3\beta^{(j)}_{3}(2)-\beta^{(j)}_{3}(0)+\beta^{(j)}_{3}(3)\right],
\end{align*}}
where one can already observe the binomial pattern that we hinted at.
\oprocend\end{example}

We will now obtain a rearrangement of \eqref{eq:nonlineards} for a general $\IPCN(M,N)$ system.

\begin{theorem}\label{thm:nonlinear}
For the $\IPCN(M,N)$ system, the approximating dynamical system~\eqref{eq:nonlineards} can be written as
\begin{align}\label{eq:nonlinear}
P_{i}(t) \approx& \sum\limits_{j=1}^{N}s_{ij}\beta^{(j)}_{M}(0)\cr
+ \sum\limits_{j=1}^{N}&\sum\limits_{n=1}^{M}\bigg[\bigg(\sum\limits_{k=0}^{n}\bigg((-1)^{n-k} {n\choose k} s_{ij}\beta^{(j)}_{M}(k)\bigg)\bigg)\bigg(\sum\limits_{\substack{(d_{1},\cdots ,d_{n})\\\in H_{n,M}}}P_{j}(t-d_{1})\cdots P_{j}(t-d_{n})\bigg)\bigg],
\end{align}
where
\[H_{n,M}:=\{(d_{1},d_{2},\cdots ,d_{n})\quad \big| \quad d_{i}\in\{1,\cdots ,M\},\quad d_{i}\neq d_{j} \quad \forall i,j\in\{1,\cdots ,n\}\}.\]
\end{theorem}
\begin{proof}
We show that \eqref{eq:nonlinear} is obtained by a rearrangement of \eqref{eq:nonlineards}. 
The R.H.S. of \eqref{eq:nonlineards} is given by:
\begin{align}\label{eq:clean}
 \sum\limits_{j=1}^{N}\bigg[\sum\limits_{B_{M}}s_{ij}\beta^{(j)}_{M}(v_{M})\prod\limits_{k=1}^{M}\bigg(a_{k}P_{j}(t-k) + (1-a_{k})(1-P_{j}(t-k))\bigg)\bigg].
 \end{align}
The constant term can be extracted from \eqref{eq:clean} by setting $\underline{a}_{j} = (0,0,\cdots ,0)$ in $B_{M}$ for $1\leq j \leq N$ and is given by
\[\sum\limits_{j=1}^{N}\frac{s_{ij}\rho_{j}}{1+M\delta_{b,j}}.\]
Now, fixing $j \in \{1,2,\cdots ,N\}$, we expand the term
\begin{equation}\label{eq:fixedj}\sum\limits_{B_{M}}s_{ij}\beta^{(j)}_{M}(v_{M})\prod\limits_{k=1}^{M}(a_{k}P_{j}(t-k) + (1-a_{k})(1-P_{j}(t-k))).
\end{equation}
Note that the order of \eqref{eq:fixedj} is $M$. 
In order to get the $n$th degree term (where $ 1\leq n \leq M $ in~\eqref{eq:fixedj}), we need to choose $n$ corresponding $P_{j}(t-k)$'s, where $ k \in \{1,2,\cdots ,M\}$ from the product 
\[
\prod\limits_{k=1}^{M}(a_{k}P_{j}(t-k) + (1-a_{k})(1-P_{j}(t-k)))
\]
and the rest $M-n$ chosen terms have to be $1$. We then look at the coefficient of the chosen $n$th order term. Note that the coefficients of the chosen $P_{j}(t-k)$'s are either~$1$ or~$-1$, depending on the tuple $\underline{a}_{j}$. Given a tuple $\underline{a}_{j}$,  there are exactly $v_{M}$, $P_{j}(t-k)$'s with coefficients $1$ and the rest $n-v_{M}$ of them have coefficient $-1$.

 Summing over all the possible coefficients of the $n$th degree term of  \eqref{eq:clean} we get
\[ \sum\limits_{k=0}^{n}(-1)^{n-k} {n\choose k} s_{ij}\beta^{(j)}_{M}(k)\sum\limits_{\substack{(d_{1},\cdots ,d_{n})\\\in H_{n,M}}}P_{j}(t-d_{1})\cdots P_{j}(t-d_{n}).\]
Finally, we can obtain the $n$th degree terms ($1\le n \le M$) for the other $N-1$ urns in exactly the same way as above.
\end{proof}

The analysis of the nonlinear dynamical systems given in~\eqref{eq:nonlinear} is clearly more intricate than the one in the case with memory one, where the evaluations were given by a linear dynamical system, namely~\eqref{eq:linear_dynamical_system}. This being said, given that the presence of nonlinearity is due to the product of probabilities, we can  use a further approximation by considering the leading linear terms. 
\begin{corollary}\label{thm:linearpart}
The linear part of the dynamical system~\eqref{eq:nonlinear} is given by
\begin{equation}\label{eq:linear,Mgeq1}P_{i}(t)\approx
\sum_{j=1}^{N}s_{ij}\beta^{(j)}_{M}(0) + \sum\limits_{j=1}^{N}\sum\limits_{k=1}^{M}s_{ij}\left(\beta^{(j)}_{M}(1)-\beta^{(j)}_{M}(0)\right)P_{j}(t-k).
\end{equation}
\end{corollary}
\begin{proof}
Setting  $n=1$ in~\eqref{eq:nonlinear} we obtain
\begin{align*}P_{i}(t) &\approx \sum\limits_{j=1}^{N}s_{ij}\beta^{(j)}_{M}(0)\\
& + \sum\limits_{j=1}^{N}\Big((-1) {1 \choose 0}s_{ij}\beta^{(j)}_{M}(0) + (-1)^{2} {1 \choose 1}s_{ij}\beta^{(j)}_{M}(k)\Big)\Big(\sum\limits_{d \in H_{1,M}}P_{j}(t-d)\Big)\\
&=\sum_{j=1}^{N}s_{ij}\beta^{(j)}_{M}(0) + \sum\limits_{j=1}^{N}\sum\limits_{k=1}^{M}s_{ij}(\beta^{(j)}_{M}(1)-\beta^{(j)}_{M}(0))P_{j}(t-k).
\end{align*}

\end{proof}
Equation~\eqref{eq:linear,Mgeq1} gives an approximate linear dynamical system for the $\IPCN(M,N)$ system. For $M\geq 1$ network of $N$ urns, we define
\[\widetilde{P}(t):= [P_{1}(t), \cdots ,P_{1}(t-M), ,P_{2}(t),\cdots ,P_{2}(t-M)\cdots ,P_{N}(t),\cdots ,P_{N}(t-M)]^{\mathbf{T}}.\]
Using \eqref{eq:nonlinear} and dropping the  nonlinear terms, we can write,
\begin{equation}\label{eq:linearsystem,Mgeq1}\widetilde{P}(t) \approx J_{N,M}\widetilde{P}(t-1) + C_{N,M}
\end{equation}
where, $J_{N,M}$ is a block matrix with $N^{2}$ blocks of size $M\times M$. 
\[J_{N,M} =\begin{bmatrix} 
J_{N,M}(1,1) & J_{N,M}(1,2) &\cdots  &  J_{N,M}(1,N)\\
J_{N,M}(2,1) & J_{N,M}(2,2) &\cdots & J_{N,M}(2,N)\\
\vdots & \ddots & \vdots& \vdots \\

J_{N,M}(N,1) & J_{N,M}(N,2) &\cdots  &J_{N,N}(N,N)
\end{bmatrix}_{NM\times NM}.\]
Here, the diagonal blocks of matrix, $J_{N,M}(i,i)$ is given by

\[\left[ 
\begin{array}{c|c} 
  \begin{array}{c c c}s_{ii}(\beta_{M}^{(i)}(1)-\beta_{M}^{(i)}(0))&\cdots& s_{ii}(\beta_{M}^{(i)}(1)-\beta_{M}^{(i)}(0))\end{array} & s_{ii}(\beta_{M}^{(i)}(1)-\beta_{M}^{(i)}(0)) \\ 
  \hline 
  \mathbf{I}_{(M-1)\times (M-1)} & \mathbf{0}_{1\times(M-1)}
  
\end{array} \right]_{M\times M}\]
where $\mathbf{I}_{(M-1)\times(M-1)}$ is the identity matrix of size $M-1$ and $\mathbf{0}_{1\times (M-1)}$ is the column vector of length $M-1$ with all entries zero. 
Similarly, the off-diagonal blocks  $J_{N,M}(i,j)$ are given by
\[\left[ 
\begin{array}{c|c} 
  \begin{array}{c c c}s_{ij}(\beta_{M}^{(j)}(1)-\beta_{M}^{(j)}(0))&\cdots& s_{ij}(\beta_{M}^{(j)}(1)-\beta_{M}^{(j)}(0))\end{array} & s_{ij}(\beta_{M}^{(j)}(1)-\beta_{M}^{(j)}(0)) \\ 
  \hline 
  \mathbf{0}_{(M-1)\times (M-1)} & \mathbf{0}_{1\times(M-1)}
\end{array} \right]_{M\times M}\]
where, $\mathbf{0}_{(M-1)\times (M-1)}$ is a matrix of size $(M-1)$ with all entries zero.
Finally, $C_{N,M}$ is a column matrix with N blocks each of size $1 \times M$ given by

\[C_{N,M}(i) = \begin{bmatrix}
\sum\limits_{j=1}^{N}s_{ij}\beta^{(j)}_{M}(0) &
0&
\cdots &
0 
\end{bmatrix}_{1\times M}^{\mathbf{T}}.
\]
 The linear dynamical system~\eqref{eq:linearsystem,Mgeq1} has a unique equilibrium which is given by $(I-J_{N,M})^{-1}C_{N,M}$. Even though we leave further studies of stability properties of the nonlinear dynamical system~\eqref{eq:nonlinear} as a future direction, it is worth pointing out that~\eqref{eq:linearsystem,Mgeq1} asymptotically converges to the unique equilibrium if and only if the spectral radius of $J_{N,M}$ is less than one. The possible dependency of this condition to the interaction matrix and urn properties is also interesting for future studies. 
 
 We next present a few simulations to assess how close these class dynamical systems are to our Markov process. 

\section{Simulation Results}\label{simulations}
\begin{figure}[ht!]
    \centering
    \includegraphics[scale=0.4]{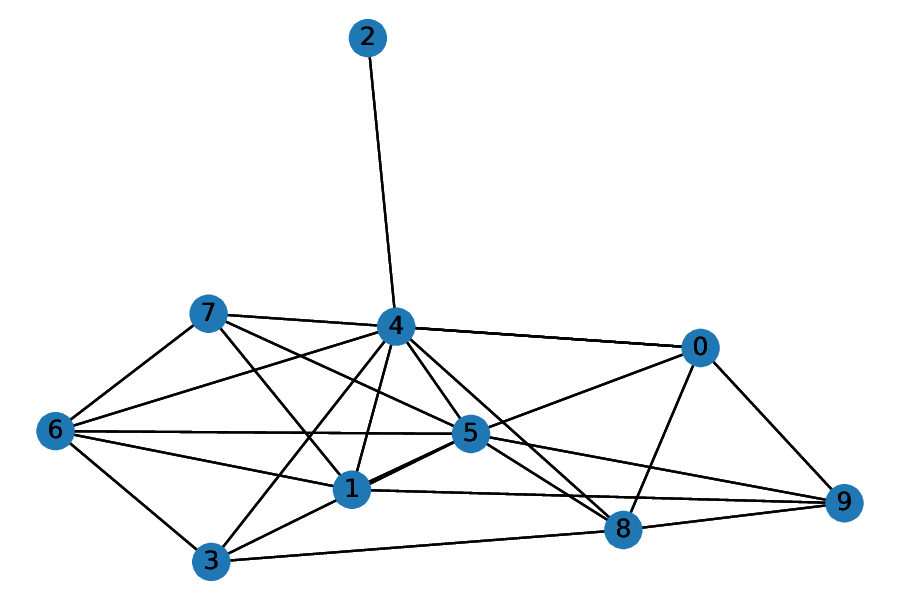}\\
    \includegraphics[scale=0.4]{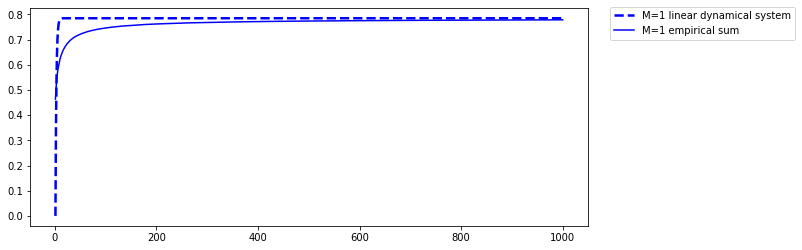}\\
\includegraphics[scale=0.4]{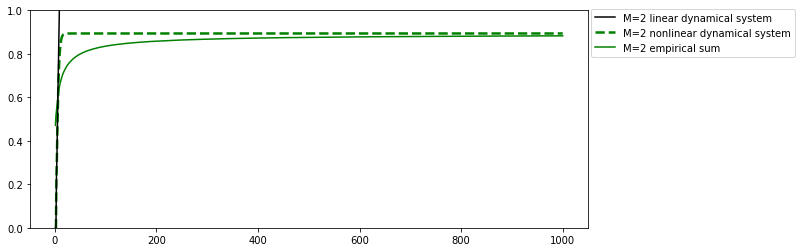}
    \includegraphics[scale=0.4]{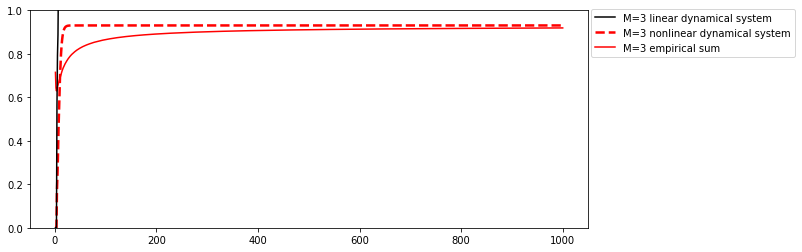}
   \caption{Infection rate curves for non-homogeneous $\IPCN(M,N)$ systems with $N=10$ nodes and memory $M=1,2,3$.
   At $t=0$, each urn has a total of $25$ balls. The number of red balls in each urn at $t=0$ is chosen randomly  between range $5$ to $23$ so that $\rho's$ lie in the range $0.2$ to $0.92$. $\Delta_{r}'s$ are chosen randomly between range $60$ to $70$ and $\Delta_{b}'s$ are randomly chosen between  range $20$ to $29$. For simplicity, we set the initial values $P_{i}(0),P_{i}(1),\cdots,P_{i}(M-1)$ all equal to zero for all urns $i$ in the network. 
   }
   \label{fig1}
\end{figure}

\begin{figure}[h!]
    \centering
    \includegraphics[scale=0.4]{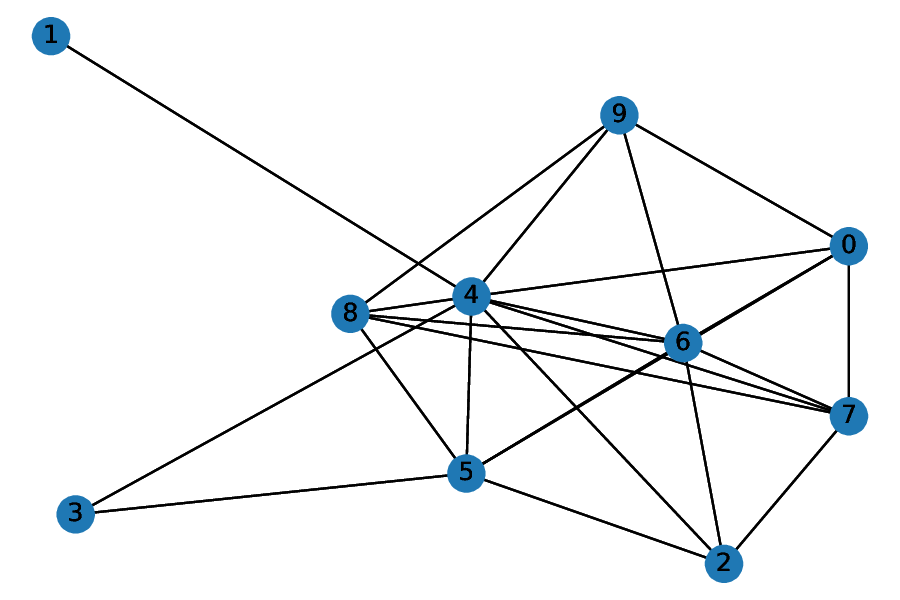}\\
    \includegraphics[scale=0.4]{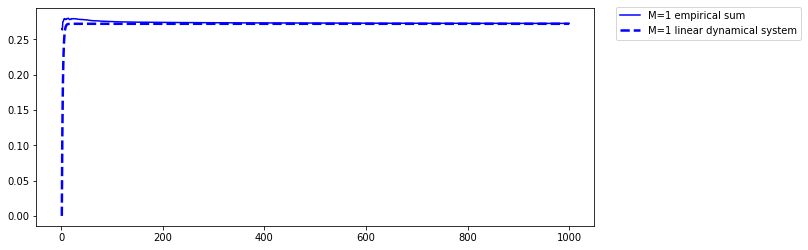}\\
    \includegraphics[scale=0.4]{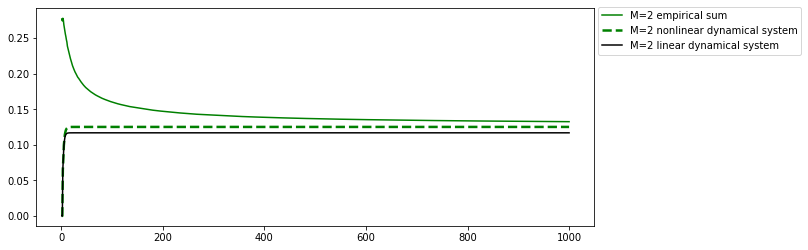}\\
    \includegraphics[scale=0.4]{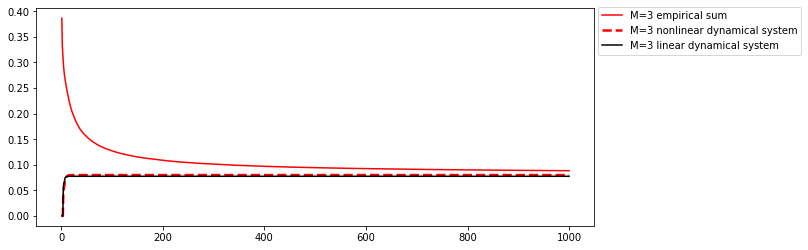}
   \caption{Infection rate curves for non-homogeneous $\IPCN(M,N)$ systems with $N=10$ nodes and memory $M=1,2,3$. At $t=0$, the total number of balls in each urn is $25$. The number of red balls in each urn at time $t=0$ are chosen randomly between the range $2$ to $17$ so that $\rho's$ lie in the range $0.08$ to $0.68$. $\Delta_{r}'s$ are chosen randomly in the range $12$ to $30$. $\Delta_{b}'s$ are chosen in the range $61$ to $80$. For simplicity, we set the initial values $P_{i}(0),P_{i}(1),\cdots,P_{i}(M-1)$ all equal to zero for all urns $i$. 
   }
\label{fig2}  
\end{figure}

\begin{figure}[h!]
    \centering
    \includegraphics[scale=0.4]{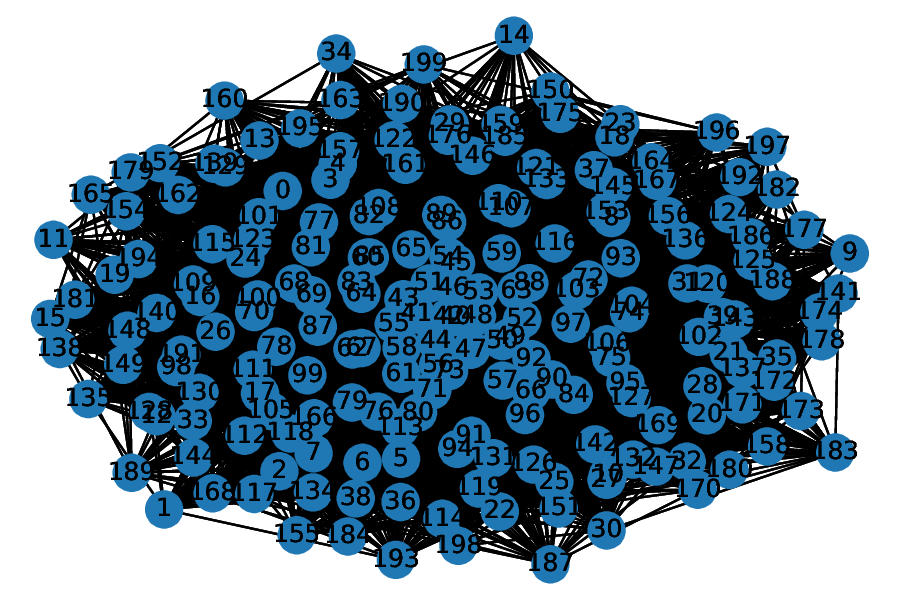}\\
    \includegraphics[scale=0.4]{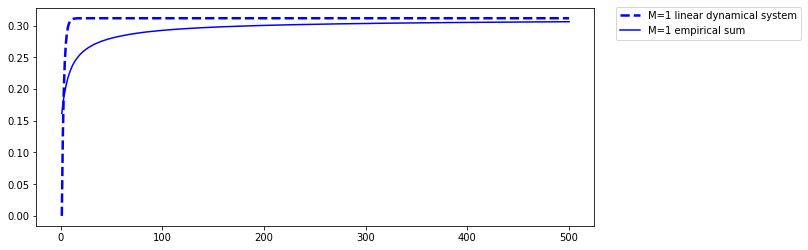}\\
    \includegraphics[scale=0.4]{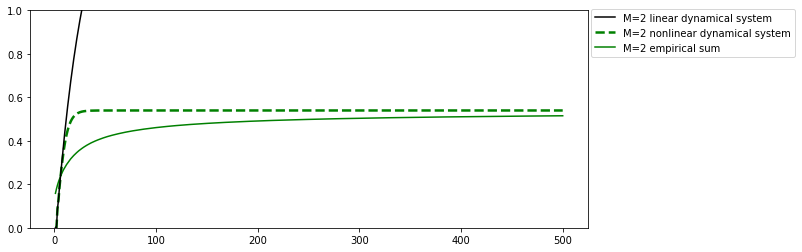}\\
    \includegraphics[scale=0.4]{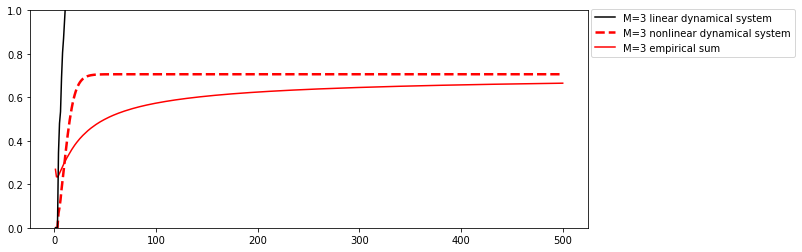}
   \caption{Infection rate curves for non-homogeneous $\IPCN(M,N)$ Barab\'{a}si-Albert \cite{RA-ALB:02} systems with $N=100$ nodes and memory $M=1,2,3$. At $t=0$, the total number of balls in each urn is $25$. The number of red balls in each urn at time $t=0$ are chosen randomly between the range $1$ to $10$ so that $\rho's$ lie in the range $0.04$ to $0.4$. $\Delta_{r}'s$ are chosen randomly in the range $40$ to $50$. $\Delta_{b}'s$ are chosen in the range $15$ to $25$. For simplicity, we set the initial values $P_{i}(0),P_{i}(1),\cdots,P_{i}(M-1)$ all equal to zero for all urns $i$. 
   }
\label{fig3}  
\end{figure}

\begin{figure}[h!]
    \centering
    \includegraphics[scale=0.3]{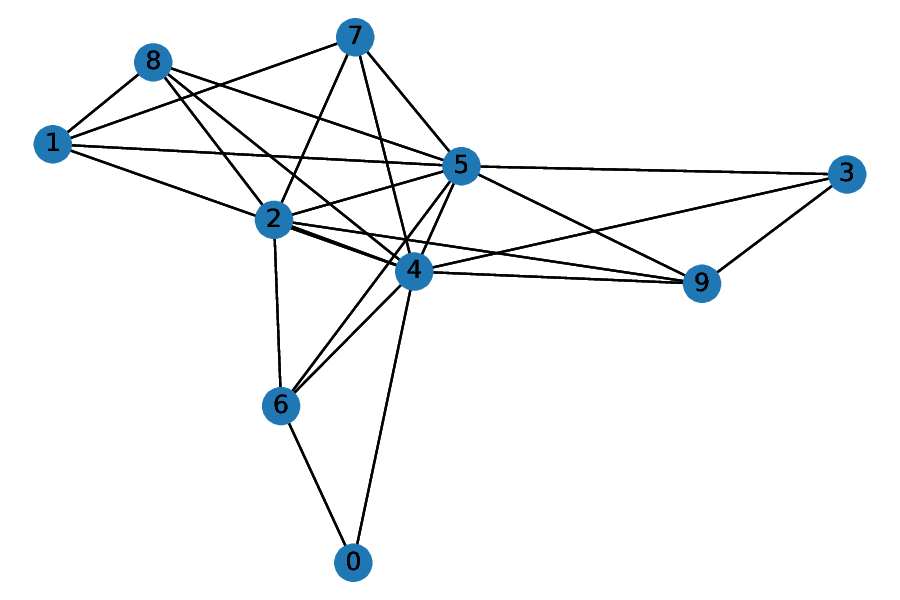}\\
    \includegraphics[scale=0.4]{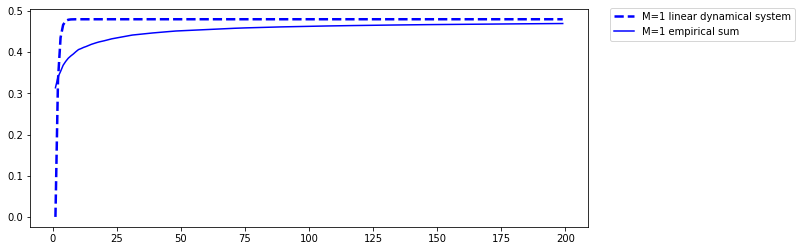}\\
    \includegraphics[scale=0.4]{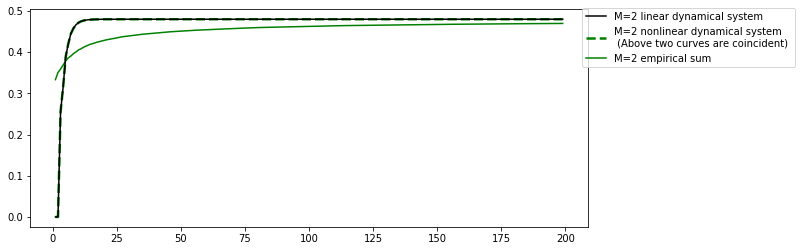}\\
    \includegraphics[scale=0.4]{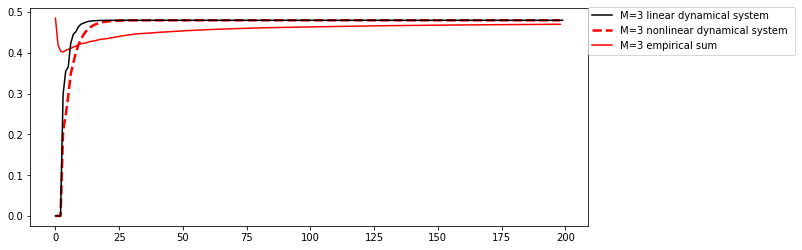}
   \caption{Infection rate curves for homogeneous $\IPCN(M,N)$ systems with $N=10$ nodes and memory $M=1,2,3$.
   We set $\rho = 0.48$,  $\delta_{r}=\delta_{b}=0.44$ for all the urns in the network. For simplicity, we set the initial values $P_{i}(0),P_{i}(1),\cdots,P_{i}(M-1)$ all equal to zero for all urns $i$. 
   }
   \label{fig4}
\end{figure}

We provide a set of simulations\footnote{For a complete list of parameters used for generating all figures, see the link: \\ \href{https://www.dropbox.com/sh/19py25reaxnfoyn/AABFdBp98J-9Jkd7zzVfTAQ9a?dl=0}{https://www.dropbox.com/sh/19py25reaxnfoyn/AABFdBp98J-9Jkd7zzVfTAQ9a?dl=0} } to illustrate our results. For this purpose, we have considered four different setups which are aimed at demonstrating the impact of memory, as well as initial urn compositions and reinforcement parameters. In particular, for the first two networks with $N=10$ (i.e.,~\cref{fig1} and~\cref{fig2}), we use $\delta_{r}$ values that are significantly larger than the $\delta_{b}$ values in~\cref{fig1} and $\delta_{b}$ values significantly larger than $\delta_{r}$ values in~\cref{fig2}. In~\cref{fig3}, we consider larger size non-homogeneous networks with $N=100$.
We simulate the $\IPCN(M,N)$ system for $M=1,2,3$ and their corresponding approximating (nonlinear) dynamical systems given by~\eqref{eq:nonlinear}. We also simulate the linear approximation~\eqref{eq:linear,Mgeq1} of the nonlinear dynamical system for each $M=2,3$. Recall that for $M=1$, the linear dynamical system in \eqref{eq:linear_dynamical_system} exactly characterizes the underlying Markov draw process.
Finally, in~\cref{fig4}, we simulate a homogeneous $\IPCN(M,N)$ system. 
Throughout, for the given $\IPCN(M,N)$ system, we plot the average empirical sum at time $ t $, which is given by
    \[
    \frac{1}{N}\sum\limits_{i=1}^{N}I_{t}(i)
    \quad \mathrm{where} \quad  
    I_{t}(i)=\frac{1}{t}\sum\limits_{n=1}^{t}Z_{i,n}.
    \]
    For each plot, the average empirical sum is computed $100$ times and the mean value is plotted against time. For the dynamical systems, we plot the average infection rate at time $t$, which is given by 
$\frac{1}{N}\sum\limits_{i=1}^{N}P_{i}(t).$


We first note from the simulations that for the network with $M=1$, the linear system in \eqref{eq:linear_dynamical_system} matches the empirical sum of the draw process, as expected since in this case the linear system is exact.

We next observe that the nonlinear dynamical system~\eqref{eq:nonlinear} is always a good approximation for the $\IPCN(M,N)$ system. 
Note that in \eqref{eq:nonlinear}, the order of the approximating nonlinear dynamical system is equal to the memory of the $\IPCN(M,N)$ system and therefore, when we drop nonlinear terms from~\eqref{eq:nonlinear} to obtain the linear approximation~\eqref{eq:linear,Mgeq1}, as we expect, the approximation gets worse. For $M>1$, we can see this worsening of linear approximation in~\cref{fig1} and~\cref{fig3}. However, in some exceptional cases, the linear approximation performs well. An example of this behavior is presented in~\cref{fig2}, where the linear approximations perform as well as the nonlinear ones. An important aspect of these simulations is that the  reinforcement parameters play a major role in determining the asymptotic value of the probability of infection. For example in~\cref{fig1}, since the $\delta_{r}$ parameters are significantly larger than the $\delta_{b}$ parameters (i.e., infection is much more likely than recovery), the asymptotic value of the plots is higher (i.e., the urns tend towards having a larger composition of red balls). Similarly in~\cref{fig2}, since the $\delta_{b}$ values are significantly larger than the $\delta_{r}$ values, the asymptotic value of the plots are lower (i.e., the urns tend towards having a larger proportion of black balls). Furthermore, the better performance of the linear system observed in \cref{fig2} relative to \cref{fig1} and \cref{fig3} is attributed to the fact that the constant term in the linear approximation (given by \eqref{eq:linear,Mgeq1}) increases when $\delta_{r}$ is increased and decreases when $\delta_{b}$ is increased. Depending on how large $\delta_{r}$ is, the probability of infection as approximated by \eqref{eq:linear,Mgeq1} can exceed $1$ and hence the linear approximation does not perform well for these cases. Whereas, no matter how large $\delta_{b}$ gets, the probability of infection never gets smaller than $0$ and hence the linear approximation performs comparatively better in this case.
Lastly, we observe from the simulations for the homogeneous $\IPCN(M,N)$ system in~\cref{fig4} that the empirical sum as well as the linear and nonlinear dynamical approximations converge to $\rho$ irrespective of the memory of the system. This phenomenon is indeed shown in Theorem~\ref{thm:homoipcn} for any homogeneous $\IPCN(M,N)$ system.

\section{Conclusions}\label{sec:conclusions}
We formulated an interacting P\'{o}lya contagion network with finite Markovian memory. We showed that for the homogeneous case, i.e., when all urns have identical initial conditions and reinforcement parameters, the underlying Markov process is irreducible and aperiodic and hence has a unique stationary distribution. We also derived the exact asymptotic marginal infection distribution. For the non-homogeneous interacting P\'{o}lya contagion network, we constructed dynamical systems to evaluate the network's infection propagation. We showed that when memory $M=1$, the probability of infection can be exactly represented by a  linear dynamical system which has a unique equilibrium point to which the solution asymptotically converges. For memory $M>1$, we used \emph{mean-field approximations} to construct approximating dynamical systems which are nonlinear in general; we obtained a linearization of this dynamical system and characterize its equilibrium. We provided simulations comparing the corresponding linear and approximating nonlinear dynamical systems with the original stochastic process. Notably, we demonstrated that the approximating nonlinear dynamical system performs well for all tested values of memory and network size.
Future work include analyzing the stability properties of the nonlinear model, studying the scaling of the approximations with the size of the network, and designing curing strategies for the proposed model with systematic comparisons with the SIS model.

\section*{Acknowledgements}
We sincerely thank Yanglei Song for providing the proof of Theorem \ref{thm:homoipcn} and improving the results in Section \ref{sec:three}.

\bibliographystyle{IEEEtran}
\bibliography{reference_arxiv}

\begin{thebibliography}{10}
\providecommand{\url}[1]{#1}
\csname url@samestyle\endcsname
\providecommand{\newblock}{\relax}
\providecommand{\bibinfo}[2]{#2}
\providecommand{\BIBentrySTDinterwordspacing}{\spaceskip=0pt\relax}
\providecommand{\BIBentryALTinterwordstretchfactor}{4}
\providecommand{\BIBentryALTinterwordspacing}{\spaceskip=\fontdimen2\font plus
\BIBentryALTinterwordstretchfactor\fontdimen3\font minus
  \fontdimen4\font\relax}
\providecommand{\BIBforeignlanguage}[2]{{%
\expandafter\ifx\csname l@#1\endcsname\relax
\typeout{** WARNING: IEEEtran.bst: No hyphenation pattern has been}%
\typeout{** loaded for the language `#1'. Using the pattern for}%
\typeout{** the default language instead.}%
\else
\language=\csname l@#1\endcsname
\fi
#2}}
\providecommand{\BIBdecl}{\relax}
\BIBdecl

\bibitem{curing}
M.~Hayhoe, F.~Alajaji, and B.~Gharesifard, ``\BIBforeignlanguage{eng}{Curing
  epidemics on networks using a {P}olya contagion model},''
  \emph{\BIBforeignlanguage{eng}{IEEE/ACM Transactions on Networking}},
  vol.~27, no.~5, pp. 2085--2097, 2019.

\bibitem{contagion}
------, ``A {P}olya contagion model for networks,'' \emph{IEEE Transactions on
  Control of Network Systems}, vol.~5, no.~4, p. 1998–2010, 12 2017.

\bibitem{consensus}
A.~{Fazeli} and A.~{Jadbabaie}, ``On consensus in a correlated model of network
  formation based on a {P}\'{o}lya urn process,'' in \emph{2011 50th IEEE
  Conference on Decision and Control and European Control Conference}, 2011,
  pp. 2341--2346.

\bibitem{AJ-AM-EM-RS:21}
A.~Jadbabaie, A.~Makur, E.~Mossel, and R.~Salhab, ``Opinion dynamics under
  social pressure,'' \emph{arXiv preprint arXiv:2104.11172v1}, 2021.

\bibitem{image}
A.~{Banerjee}, P.~{Burlina}, and F.~{Alajaji}, ``Image segmentation and
  labeling using the {P}olya urn model,'' \emph{IEEE Transactions on Image
  Processing}, vol.~8, no.~9, pp. 1243--1253, 1999.

\bibitem{internet}
N.~Berger, C.~Borgs, J.~Chayes, and A.~Saberi, ``\BIBforeignlanguage{eng}{On
  the spread of viruses on the internet},'' vol.~5.\hskip 1em plus 0.5em minus
  0.4em\relax Proceedings of the ACM-SIAM Symposium on Discrete Algorithms,
  2005, pp. 301--310.

\bibitem{social}
B.~Skyrms and R.~Pemantle, ``\BIBforeignlanguage{eng}{A dynamic model of social
  network formation},'' \emph{\BIBforeignlanguage{eng}{Proceedings of the
  National Academy of Sciences}}, vol.~97, no.~16, pp. 9340--9346, 2000.

\bibitem{neeraja}
N.~Sahasrabudhe, ``Synchronization and fluctuation theorems for interacting
  {F}riedman urns,'' \emph{J. Appl. Prob.}, vol.~53, pp. 1221--1239, 2016.

\bibitem{neeraja1}
G.~Kaur and S.~Neeraja, ``Interacting urns on a finite directed graph,'' 2019.

\bibitem{collevecchio2013preferential}
A.~Collevecchio, C.~Cotar, and M.~LiCalzi, ``On a preferential attachment and
  generalized {P}{\'o}lya’s urn model,'' \emph{The Annals of Applied
  Probability}, vol.~23, no.~3, pp. 1219--1253, 2013.

\bibitem{pa}
N.~Berger, C.~Borgs, J.~T. Chayes, and A.~Saberi, ``Asymptotic behavior and
  distributional limits of preferential attachment graphs,'' \emph{Annals of
  Probability}, vol.~42, no.~1, pp. 1--40, 01 2014.

\bibitem{greg}
G.~Harrington, F.~Alajaji, and B.~Gharesifard, ``Initialization and curing
  policies for {P}\'{o}lya contagion networks,'' \emph{SIAM Journal on Control
  and Optimization, to appear}.

\bibitem{synchronization}
P.~Pra, P.~Y. Louis, and I.~G. Minelli, ``Synchronazation via interacting
  reinforcement,'' \emph{Journal of Applied Probability}, vol.~51, pp.
  556--568, 03 2016.

\bibitem{ben}
M.~Benaïm, I.~Benjamini, J.~Chen, and Y.~Lima, ``\BIBforeignlanguage{eng}{A
  generalized {P}ólya's urn with graph based interactions},''
  \emph{\BIBforeignlanguage{eng}{Random Structures \& Algorithms}}, vol.~46,
  no.~4, pp. 614--634, 2015.

\bibitem{finite}
F.~Alajaji and T.~Fuja, ``A communication channel modeled on contagion,''
  \emph{IEEE Transactions on Information Theory}, vol.~40, no.~6, pp.
  2035--2041, 1994.

\bibitem{hassibi}
N.~A. Ruhi, T.~Christos, and B.~Hassibi, ``\BIBforeignlanguage{eng}{Improved
  bounds on the epidemic threshold of exact {SIS} models on complex
  networks},'' \emph{\BIBforeignlanguage{eng}{2016 IEEE 55th Conference on
  Decision and Control (CDC)}}, pp. 3560--3565, 2016.

\bibitem{wang}
{Yang Wang}, D.~{Chakrabarti}, {Chenxi Wang}, and C.~{Faloutsos}, ``Epidemic
  spreading in real networks: {A}n eigenvalue viewpoint,'' in \emph{Proceedings
  22nd International Symposium on Reliable Distributed Systems}, 2003, pp.
  25--34.

\bibitem{virus}
P.~Van~Mieghem, J.~Omic, and R.~Kooij, ``\BIBforeignlanguage{eng}{Virus spread
  in networks},'' \emph{\BIBforeignlanguage{eng}{IEEE/ACM Transactions on
  Networking}}, vol.~17, no.~1, pp. 1--14, 2009.

\bibitem{pare2017epidemic}
P.~E. Par{\'e}, C.~L. Beck, and A.~Nedi{\'c}, ``Epidemic processes over
  time-varying networks,'' \emph{IEEE Transactions on Control of Network
  Systems}, vol.~5, no.~3, pp. 1322--1334, 2017.

\bibitem{mei2017dynamics}
W.~Mei, S.~Mohagheghi, S.~Zampieri, and F.~Bullo, ``On the dynamics of
  deterministic epidemic propagation over networks,'' \emph{Annual Reviews in
  Control}, vol.~44, pp. 116--128, 2017.

\bibitem{nowzari2016analysis}
C.~Nowzari, V.~M. Preciado, and G.~J. Pappas, ``Analysis and control of
  epidemics: A survey of spreading processes on complex networks,'' \emph{IEEE
  Control Systems Magazine}, vol.~36, no.~1, pp. 26--46, 2016.

\bibitem{andersson2012stochastic}
H.~Andersson and T.~Britton, \emph{Stochastic {E}pidemic {M}odels and their
  {S}tatistical {A}nalysis}.\hskip 1em plus 0.5em minus 0.4em\relax Springer
  Science \& Business Media, 2012, vol. 151.

\bibitem{irene_markov}
G.~Aletti and I.~Crimaldi, ``The rescaled {P}\'{o}lya urn: {L}ocal
  reinforcement and {C}hi-squared goodness of fit test,'' 2019.

\bibitem{polya-lundberg}
D.~Pfeifer, \emph{P\'{o}lya-{L}undberg process}.\hskip 1em plus 0.5em minus
  0.4em\relax Encyclopedia of Statistical Sciences, vol.~7, Wiley, New York,
  1986.

\bibitem{corona}
N.~R. Barraza, G.~Pena, and V.~Moreno, ``A non-homogeneous {M}arkov early
  epidemic growth dynamics model. application to the {SARS}-{C}o{V}-2
  pandemic,'' \emph{Chaos, Solitons \& Fractals}, p. 110297, 2020.

\bibitem{FCF:21}
F.~C. Fabiani, ``Asymptotic incidence rate estimation of
  {S}{A}{R}{S}-{C}{O}{V}{I}{D}-19 via a {P}\'{o}lya process scheme: {A}
  comparative analysis in {I}taly and {E}uropean countries,'' \emph{arXiv
  preprint arXiv:2010.00463}, 2020.

\bibitem{friedman}
D.~A. Freedman, ``Bernard {F}riedman's urn,'' \emph{The Annals of Mathematical
  Statistics}, vol.~36, no.~3, pp. 956--970, 1965.

\bibitem{martingale}
R.~Gouet, ``Martingale functional central limit theorems for a generalized
  {P}\'{o}lya urn,'' \emph{The Annals of Probability}, vol.~21, no.~3, pp.
  1624--1639, 07 1993.

\bibitem{sa}
V.~S. Borkar, \emph{Stochastic Approximation: A Dynamical Systems Viewpoint},
  ser. Texts and readings in Mathematics.\hskip 1em plus 0.5em minus
  0.4em\relax Cambridge University Press, 2008.

\bibitem{krishanu}
U.~Gangopadhyay and K.~Maulik, ``Stochastic approximation with random step
  sizes and urn models with random replacement matrices,'' \emph{Annals of
  Applied Probability}, vol.~29, no.~4, pp. 2033--2066, 2019.

\bibitem{sophie}
S.~Laurelle and G.~Pages, ``Randomized urn models revisited using stochastic
  approximation,'' \emph{The Annals of Applied Probability}, vol.~23, no.~4,
  pp. 1409--1436, 06 2013.

\bibitem{bp}
K.~B. Athreya and P.~E. Ney, \emph{{B}ranching {P}rocesses}, ser. Dover Books
  on Mathematics.\hskip 1em plus 0.5em minus 0.4em\relax Dover Publications,
  2004.

\bibitem{embedding}
K.~B. Athereya and S.~Karlin, ``Embedding of urn schemes into continuous time
  {M}arkov branching processes and related limit theorems,'' \emph{The Annals
  of Mathematical Statistics}, vol.~39, no.~6, pp. 1801--1817, 1968.

\bibitem{branching(1)}
S.~Janson, ``Functional limit theorems for multitype branching processes and
  generalized {P}ólya urns,'' \emph{Stochastic Processes and their
  Applications}, vol. 110, no.~2, pp. 177 -- 245, 2004.

\bibitem{survey}
R.~Pemantle, ``A survey of random processes with reinforcement,'' \emph{Probab.
  Surveys}, vol.~4, pp. 1--79, 2007.

\bibitem{irene}
G.~Aletti, I.~Crimaldi, and A.~Ghiglietti, ``Synchronization of reinforced
  stochastic processes with a network based interaction,'' \emph{The Annals of
  Applied Probability}, vol.~27, no.~6, pp. 3787--3844, 2017.

\bibitem{spectrum}
S.~H. Strogatz, \emph{\BIBforeignlanguage{eng}{Nonlinear {D}ynamics and
  {C}haos: {W}ith {A}pplications to {P}hysics, {B}iology, {C}hemistry, and
  {E}ngineering}}, 2nd~ed.\hskip 1em plus 0.5em minus 0.4em\relax CRC Press,
  2018.

\bibitem{RA-ALB:02}
R.~Albert and A.~L. Barab{\'a}si, ``Statistical mechanics of complex
  networks,'' \emph{Rev.\ Modern Phys.}, vol.~74, no.~1, pp. 47--97, 2002.

\end{thebibliography}


\clearpage

\end{spacing}
\end{document}